\documentclass[a4paper,11pt]{amsart}
\usepackage{graphicx}
\usepackage{mathptmx}
\usepackage{amsmath}
\usepackage{amssymb}
\usepackage{enumitem}
\usepackage{xcolor}

\newmuskip\pFqmuskip

\newcommand*\pFq[6][8]{%
  \begingroup % only local assignments
  \pFqmuskip=#1mu\relax
  \mathcode`=\string"8000
  \begingroup\lccode`\~=`\,
  \lowercase{\endgroup\let~}\pFqcomma
  F^{#2}_{#3}{\left(\genfrac..{0pt}{}{#4}{#5}\bigg|#6\right)}%
  \endgroup
}
\newcommand{\pFqcomma}{\mskip\pFqmuskip}

\newtheorem{theorem}{Theorem}
\newtheorem{lemma}[theorem]{Lemma}
\newtheorem{corollary}[theorem]{Corollary}
\newtheorem{proposition}[theorem]{Proposition}

\begin{document}

\title[A note on truncated degenerate Bell polynomials]{A note on truncated degenerate Bell polynomials}

\author{Taekyun  Kim}
\address{Department of Mathematics, Kwangwoon University, Seoul 139-701, Republic of Korea}
\email{tkkim@kw.ac.kr}

\author{DAE SAN KIM}
\address{Department of Mathematics, Sogang University, Seoul 121-742, Republic of Korea}
\email{dskim@sogang.ac.kr}

\subjclass[2010]{11B73; 11B83; 65C50}
\keywords{truncated degenerate Bell polynomials; truncated degenerate modified Bell polynomials; beta random variable}

\begin{abstract}
The aim of this paper is to introduce truncated degenerate Bell polynomials and numbers and to investigate some of  their properties. In more detail, we obtain explicit expressions, identities involving other special polynomials, integral representations, Dobinski-like formula and expressions of the generating function in terms of differential operators and linear incomplete gamma function. In addition, we introduce  truncated degenerate modified Bell polynomials and numbers and get similar results for those polynomials.  As an application of our results, we show that the truncated degenerate Bell numbers can be expressed as a finite sum involving moments of a beta random variable with certain parameters.
\end{abstract}

\maketitle

\section{Introduction}
Carlitz is the first one who initiated a study of degenerate versions of some special polynomials and numbers, namely degenerate Bernoulli polynomials and numbers and degenerate Euler polynomials and numbers. In recent years, extensive researches have been done for various degenerate versions of some special polynomials and numbers and have yielded many interesting arithmetical and combinatorial results.
These include the degenerate Stirling numbers of the first and second kinds, degenerate central factorial numbers of the second kind, degenerate Bernoulli numbers of the second kind, degenerate Bernstein polynomials, degenerate Bell numbers and polynomials, degenerate central Bell numbers and polynomials, degenerate complete Bell polynomials and numbers, and so on. \par
Truncated polynomials have been shown to play an important role in various areas. However, not much are known for the properties of these polynomials. The following are some of the old results related to such polynomials. The minimum variance unbiased estimation is discussed in [4] for the zero class truncated bivariate Poisson and logarithmic series distributions, the maximum likelihood estimation of the Poisson parameter $\lambda$ is concerned in [5] when the zero class has been truncated and a new family of Hermite polynomials is constructed in [7] by using the truncated exponential with applications to flattened beams in optics. More recently,  degenerate exponential  truncated polynomials and numbers are studied in [9], the degenerate zero-truncated Poisson random variables are introduced in [14], the truncated-exponential-based Apostol-type polynomials are investigated in [23] and the truncated exponential-based Mittag-Leffler polynomials are considered in [25]. Further, in [12] an umbral calculus approach is given for Bernoulli-Pad\'e polynomials of fixed order, which include the truncated Bernoulli polynomials as a special case and whose generating function is based on the Pad\'e approximant of the exponential function. \par
The aim of this paper is to introduce truncated degenerate Bell polynomials and numbers (see \eqref{13}) and to explore their various properties. In more detail, for the truncated degenerate Bell polynomials and numbers we obtain explicit expressions, identities involving other special polynomials, integral representations, Dobinski-like formula and expressions of the generating function in terms of differential operators and  linear incomplete gamma function. In addition, we introduce  truncated degenerate modified Bell polynomials (see \eqref{41}) and numbers and get similar results for those polynomials.  Finally, as an application of our results, we show that, if $X$ is the beta random variable with parameters $1,p$, then the truncated degenerate Bell numbers can be expressed as a finite sum involving moments of $X$. In the rest of this section, we recall some necessary facts that are needed throughout this paper. \par
For any $\lambda\in\mathbb{R}$, the degenerate exponential functions are defined by
\begin{equation}
e_{\lambda}^{x}(t)=\sum_{n=0}^{\infty}(x)_{n,\lambda}\frac{t^{n}}{n!},\quad e_{\lambda}(t)=e_{\lambda}^{1}(t),\quad(\mathrm{see}\ [16,17]),\label{1}
\end{equation}
where $(x)_{0,\lambda}=1$, $(x)_{n,\lambda}=x(x-\lambda)\cdots(x-(n-1)\lambda)$, $(n\ge 1)$.\par  Note that $\displaystyle\lim_{\lambda\rightarrow 0}e_{\lambda}^{x}(t)=e^{xt}\displaystyle$. The Stirling numbers of the first kind $S_{1}(n,k)$ appear as the coefficients in the expansion
\begin{equation}
(x)_{n}=\sum_{k=0}^{n}S_{1}(n,k)x^{k},\quad(n\ge 0),\quad(\mathrm{see}\ [6,20]),\label{2}
\end{equation}
where $(x)_{0}=1$, $(x)_{n}=x(x-1)\cdots(x-n+1)$, $(n\ge 1)$. \par
As the inversion formula of \eqref{2}, the Stirling numbers of the second kind $S_{2}(n,k)$ are defined as the coefficients in the expansion
\begin{equation}
x^{n}=\sum_{k=0}^{n}S_{2}(n,k)(x)_{k},\quad(n\ge 0),\quad(\mathrm{see}\ [6,20]).\label{3}
\end{equation}
The Stirling number of the second kind $S_{2}(n,k)$ counts the number of ways of partitioning a set of $n$ elements into $k$ non-empty subsets. The number of all partitions of a set of $n$ elements is the Bell number $\mathrm{Bel}_{n},\ (n\ge 0)$. Thus we note that
\begin{equation}
\mathrm{Bel}_{n}=\sum_{k=0}^{n}S_{2}(n,k),\quad(n\ge 0),\quad(\mathrm{see}\ [6,20]).\label{4}
\end{equation}
Further, the Bell polynomials are given by
\begin{equation}
\mathrm{Bel}_{n}(x)=\sum_{k=0}^{n}S_2(n,k)x^k, \label{4-1}
\end{equation}
with the generating function
\begin{equation}
e^{x(e^t-1)}=\sum_{n=0}^{\infty}\mathrm{Bel}_n(x)\frac{t^n}{n!}, \quad(\mathrm{see}\ [2,3,8]). \label{4-2}
\end{equation} \par
Many researchers have studied Bell numbers in connection with various areas (see [1,2,3,8,15,18]). Recently, the degenerate Stirling numbers of the first kind are defined as the coefficients in the expansion
\begin{equation}
(x)_{n}=\sum_{k=0}^{n}S_{1,\lambda}(n,k)(x)_{k,\lambda},\quad(n\ge 0),\quad(\mathrm{see}\ [1,11]).\label{5}
\end{equation}
As the inversion formula of \eqref{5}, the degenerate Stirling numbers of the second kind appear as the coefficients in the expansion
\begin{equation}
(x)_{n,\lambda}=\sum_{k=0}^{n}S_{2,\lambda}(n,k)(x)_{k},\quad(n\ge 0),\quad(\mathrm{see}\ [10]).\label{6}
\end{equation}
Let $\log_{\lambda}t$ be the compositional inverse of $e_{\lambda}(t)$. Then we note that the generating function of the degenerate Stirling numbers of the first kind is given by
\begin{equation}
\frac{1}{k!}\big(\log_{\lambda}(1+t)\big)^{k}=\sum_{n=k}^{\infty}S_{1,\lambda}(n,k)\frac{t^{n}}{n!},\quad(\mathrm{see}\ [10,13]).\label{7}
\end{equation}
It is well known that the degenerate Bernoulli polynomials of order $r\ (\in\mathbb{N})$ are given by
\begin{equation}
\bigg(\frac{t}{e_{\lambda}(t)-1}\bigg)^{r}e_{\lambda}^{x}(t)=\sum_{n=0}^{\infty}\beta_{n,\lambda}^{(r)}(x)\frac{t^{n}}{n!},\quad(\mathrm{see}\ [10,13]).\label{8}
\end{equation}
In [18], the degenerate Bell polynomials are defined by
\begin{equation}
e^{x(e_{\lambda}(t)-1)}=\sum_{n=0}^{\infty}\mathrm{Bel}_{n,\lambda}(x)\frac{t^{n}}{n!},\quad(\mathrm{see}\ [15,19]). \label{9}
\end{equation}
When $x=1$, $\mathrm{Bel}_{n,\lambda}=\mathrm{Bel}_{n,\lambda}(1)$ are called the degenerate Bell numbers.  Note that
\begin{displaymath}
	\lim_{\lambda\rightarrow 0}\mathrm{Bel}
_{n,\lambda}=\mathrm{Bel}_{n},\quad (n\ge 0).
\end{displaymath}
As in [17,18], we note that
\begin{equation}
\mathrm{Bel}_{n,\lambda}(x)=\sum_{k=0}^{n}S_{2,\lambda}(n,k)x^k,\quad(n\ge 0).\label{10}
\end{equation}
As is well known, the beta function is given by the integral
\begin{equation}
B(\alpha,\beta)=\int_{0}^{1}t^{\alpha-1}(1-t)^{\beta-1}dt,\quad(\mathrm{Re}\,\alpha,\,\mathrm{Re}\,\beta>0),\quad(\mathrm{see}\ [24]).\label{11}
\end{equation}
The beta function in \eqref{11} is equal to
\begin{equation}
B(\alpha,\beta)=\frac{\Gamma(\alpha)\Gamma(\beta)}{\Gamma(\alpha+\beta)},\quad(\mathrm{see}\ [24]),\label{12}
\end{equation}
where $\Gamma(\alpha)$ is the gamma function with the property $\Gamma(\alpha+1)=\alpha\Gamma(\alpha)$. \par

\section{Truncated degenerate Bell numbers and polynomials}
Truncating and then normalizing the degenerate Bell polynomials in \eqref{9}, we define {\it{truncated degenerate Bell polynomials}} $\mathrm{Bel}_{n,\lambda}^{(p)}(x)$ by
\begin{equation}
\sum_{n=0}^{\infty}\mathrm{Bel}_{n,\lambda}^{(p)}(x)\frac{t^{n}}{n!}=\frac{p!}{x^{p}(e_{\lambda}(t)-1)^{p}}\bigg(e^{x(e_{\lambda}(t)-1)}-\sum_{k=0}^{p-1}\frac{x^{k}(e_{\lambda}(t)-1)^{k}}{k!}\bigg),\label{13}
\end{equation}
where $p$ is a nonnegative integer. In addition, $\mathrm{Bel}_{n,\lambda}^{(p)}=\mathrm{Bel}_{n,\lambda}^{(p)}(1)$ are called {\it{truncated degenerate Bell numbers}}.
\par
From \eqref{13}, it is immediate to see that
\begin{equation}\label{14}
\sum_{n=0}^{\infty}\mathrm{Bel}_{n,\lambda}^{(p)}(x)\frac{t^{n}}{n!}= p!\sum_{k=0}^{\infty}\frac{x^{k}\big(e_{\lambda}(t)-1\big)^{k}}{(k+p)!}.
\end{equation}
From \eqref{14}, we observe that
\begin{align}
\sum_{n=0}^{\infty}\mathrm{Bel}_{n,\lambda}^{(p)}(x)\frac{t^{n}}{n!}\
&=\ \sum_{k=0}^{\infty}\frac{p!k!}{(k+p)!}x^{k}\frac{1}{k!}\big(e_{\lambda}(t)-1\big)^{k}\label{15} \\
&=\ \sum_{k=0}^{\infty}\frac{x^{k}}{\binom{k+p}{k}}\sum_{n=k}^{\infty}S_{2,\lambda}(n,k)\frac{t^{n}}{n!}\nonumber \\
&=\ \sum_{n=0}^{\infty}\Bigg(\sum_{k=0}^{n}\frac{S_{2,\lambda}(n,k)}{\binom{k+p}{k}}x^{k}\Bigg)\frac{t^{n}}{n!}.\nonumber	
\end{align}
Therefore, by comparing the coefficients on both sides \eqref{15}, we obtain the following theorem.
\begin{theorem}
For $n,p\ge 0$, we have
\begin{displaymath}
\mathrm{Bel}_{n,\lambda}^{(p)}(x) = \sum_{k=0}^{n}\frac{S_{2,\lambda}(n,k)}{\binom{k+p}{k}}x^{k}.
\end{displaymath}
\end{theorem}
When $p=0$, we have $\mathrm{Bel}_{n,\lambda}^{(0)}(x)=\mathrm{Bel}_{n,\lambda}(x),\ (n\ge 0)$.
Let us take $p=1$. Then, by \eqref{15}, we get

\begin{align}
\sum_{n=0}^{\infty}\mathrm{Bel}_{n,\lambda}^{(1)}(x)\frac{t^{n}}{n!}\
&=\ \frac{1}{x}\frac{t}{e_{\lambda}(t)-1}\frac{1}{t}\big(e^{x(e_{\lambda}(t)-1)}-1\big)  \label{16}   \\
&=\ \frac{1}{x}\sum_{l=0}^{\infty}\beta_{l,\lambda}\frac{t^{l}}{l!}\frac{1}{t}\sum_{m=1}^{\infty}\mathrm{Bel}_{m,\lambda}(x)\frac{t^{m}}{m!}\nonumber \\
&=\frac{1}{x}\sum_{n=1}^{\infty}\sum_{m=1}^{n}\frac{1}{t}\binom{n}{m}\beta_{n-m,\lambda}Bel_{m,\lambda}(x)\frac{t^n}{n!} \nonumber \\
&=\frac{1}{x}\sum_{n=0}^{\infty}\sum_{m=0}^{n}\frac{1}{t}\binom{n+1}{m+1}\beta_{n-m,\lambda}Bel_{m+1,\lambda}(x)\frac{t^{n+1}}{(n+1)!} \nonumber \\
&=\ \sum_{n=0}^{\infty}\bigg(\frac{1}{x}\sum_{m=0}^{n}\frac{1}{m+1}\binom{n}{m}\beta_{n-m,\lambda}\mathrm{Bel}_{m+1,\lambda}(x)\bigg)\frac{t^{n}}{n!},\nonumber
\end{align}
where $\beta_{n,\lambda}=\beta_{n,\lambda}^{(1)}(0)$ are called the Carlitz's degenerate Bernoulli numbers. \par
Therefore, by \eqref{16}, we obtain the following theorem.
\begin{theorem}
For $n\ge 0$, we have
\begin{displaymath}
x\mathrm{Bel}_{n,\lambda}^{(1)}(x)= \sum_{m=0}^{n}\frac{1}{m+1}\binom{n}{m}\beta_{n-m,\lambda}\mathrm{Bel}_{m+1,\lambda}(x).
\end{displaymath}
In the case of $x=1$, we obtain
\begin{displaymath}
\mathrm{Bel}_{n,\lambda}^{(1)}=\sum_{m=0}^{n}\frac{1}{m+1}\binom{n}{m}\beta_{n-m,\lambda}\mathrm{Bel}_{m+1,\lambda}.
\end{displaymath}
\end{theorem}
We observe from \eqref{14} that
\begin{align}
p\int_{0}^{1}e^{x(e_{\lambda}(t)-1)}(1-x)^{p-1}dx \ &=\ p\sum_{n=0}^{\infty}\frac{(e_{\lambda}(t)-1)^{n}}{n!}\int_{0}^{1}x^{n}(1-x)^{p-1}dx \label{17} \\
&=\ p\sum_{n=0}^{\infty}\frac{(e_{\lambda}(t)-1)^{n}}{n!}B(n+1,p) \nonumber \\
&=\ p!\sum_{n=0}^{\infty}\frac{(e_{\lambda}(t)-1)^{n}}{(n+p)!}\nonumber\\
&= \sum_{n=0}^{\infty}\mathrm{Bel}_{n,\lambda}^{(p)}\frac{t^{n}}{n!}. \nonumber
\end{align}
On the other hand, we also have
\begin{equation}
p\int_{0}^{1}e^{x(e_{\lambda}(t)-1)}(1-x)^{p-1}dx=\sum_{n=0}^{\infty}p\int_{0}^{1}\mathrm{Bel}_{n,\lambda}(x)(1-x)^{p-1}dx\frac{t^{n}}{n!}. \label{18}
\end{equation}
Therefore, by \eqref{17} and \eqref{18}, we obtain the following proposition.
\begin{proposition}
For $p,n\ge 0$, we have
\begin{displaymath}
p\int_{0}^{1}\mathrm{Bel}_{n,\lambda}(x)(1-x)^{p-1}dx=\mathrm{Bel}_{n,\lambda}^{(p)}.
\end{displaymath}
\end{proposition}
By \eqref{10} and Proposition 3, we get
\begin{align}
\mathrm{Bel}_{n,\lambda}^{(p)}\ &=\ p\sum_{k=0}^{n}S_{2,\lambda}(n,k)\int_{0}^{1}x^{k}(1-x)^{p-1}dx \label{19}\\
&= p\sum_{k=0}^{n}S_{2,\lambda}(n,k)B(k+1,p)= \sum_{k=0}^{n}\frac{S_{2,\lambda}(n,k)}{\binom{k+p}{k}},\nonumber
\end{align}
which also follows from Theorem 1 with $x=1$.
From \eqref{17}, we note that
\begin{align}
\sum_{n=0}^{\infty}\mathrm{Bel}_{n,\lambda}^{(p)}\frac{t^{n}}{n!}
&= p\int_{0}^{1}e^{x(e_{\lambda}(t)-1)}(1-x)^{p-1}dx \label{20} \\
&= p\int_{0}^{1}e^{xe_{\lambda}(t)}e^{-x}(1-x)^{p-1}dx \nonumber \\
&= p\sum_{k=0}^{\infty}\frac{1}{k!}e_{\lambda}^{k}(t)\int_{0}^{1}e^{-x}x^{k}(1-x)^{p-1}dx  \nonumber \\
&= \sum_{n=0}^{\infty}\bigg(p\sum_{k=0}^{\infty}\frac{(k)_{n,\lambda}}{k!}\int_{0}^{1}e^{-x}(1-x)^{p-1}x^{k}dx\bigg)\frac{t^{n}}{n!}. \nonumber
\end{align}
For $n\ge 0$, by \eqref{20}, we get
\begin{align}
\mathrm{Bel}_{n,\lambda}^{(p)}\ &=\ p\sum_{k=0}^{\infty}\frac{(k)_{n,\lambda}}{k!}\int_{0}^{1}e^{-x}(1-x)^{p-1}x^{k}dx \label{21} \\
&=\ p\sum_{k=0}^{\infty}\frac{(k)_{n,\lambda}}{k!}\sum_{l=0}^{\infty}\frac{(-1)^{l}}{l!}\int_{0}^{1}(1-x)^{p-1}x^{l+k}dx \nonumber \\
&=\ p\sum_{k=0}^{\infty}\sum_{l=0}^{\infty}\frac{(k)_{n,\lambda}}{k!}\frac{(-1)^{l}}{l!}B(p,l+k+1)\nonumber\\
&=\ \sum_{k=0}^{\infty}\sum_{l=0}^{\infty}(-1)^{l}\frac{\binom{k+l}{l}}{\binom{k+l+p}{p}}\frac{(k)_{n,\lambda}}{(k+l)!}.\nonumber	
\end{align}
Therefore, by \eqref{21}, we obtain the following Dobinski-like formula.
\begin{theorem}[Dobinski-like formula]
For $n\ge 0$, we have
\begin{displaymath}
\mathrm{Bel}_{n,\lambda}^{(p)}=\sum_{k=0}^{\infty}\sum_{l=0}^{\infty}(-1)^{l}\frac{\binom{k+l}{l}}{\binom{k+l+p}{p}}\frac{(k)_{n,\lambda}}{(k+l)!}.
\end{displaymath}
\end{theorem}
Note that
\begin{displaymath}
\lim_{\lambda\rightarrow 0}\mathrm{Bel}_{n,\lambda}^{(0)}=\sum_{k=0}^{\infty}\sum_{l=0}^{\infty}\frac{(-1)^{l}}{l!k!}k^{n}=\frac{1}{e}\sum_{k=0}^{\infty}\frac{k^{n}}{k!}=\mathrm{Bel}_{n},\quad(n\ge 0).
\end{displaymath}
For $n\ge 0$, by Proposition 3, we get
\begin{align}
\mathrm{Bel}_{n,\lambda}^{(p)}\ &=\ p\int_{0}^{1}\mathrm{Bel}_{n,\lambda}(x)(1-x)^{p-1}dx	\ =\ p\sum_{k=0}^{n}S_{2,\lambda}(n,k)\int_{0}^{1}x^{k}(1-x)^{p-1}dx \label{22}\\
&=\ p\sum_{k=0}^{n}S_{2,\lambda}(n,k)\sum_{m=0}^{p-1}(-1)^{m}\binom{p-1}{m}\int_{0}^{1}x^{k+m}dx \nonumber \\
&=\ \sum_{k=0}^{n}\sum_{m=0}^{p-1}\frac{(m+1)\binom{p}{m+1}(-1)^{m}S_{2,\lambda}(n,k)}{k+m+1}.\nonumber
\end{align}
By \eqref{17}, we see that
\begin{align}
\sum_{n=0}^{\infty}\mathrm{Bel}_{n,\lambda}^{(p)}\frac{t^{n}}{n!}\ &=\ p\int_{0}^{1}e^{(e_{\lambda}(t)-1)x}(1-x)^{p-1}dx  \label{23} \\
&=\ p\int_{0}^{1}e^{(e_{\lambda}(t)-1)(1-x)}x^{p-1}dx  \nonumber \\
&=\ p\frac{e^{e_{\lambda}(t)-1}}{\big(e_{\lambda}(t)-1\big)^{p}}\int_{0}^{e_{\lambda}(t)-1}e^{-y}y^{p-1}dy \nonumber \\
&=\ p\frac{e^{e_{\lambda}(t)-1}}{\big(e_{\lambda}(t)-1\big)^{p}}d\big(p,e_{\lambda}(t)-1\big),\nonumber
\end{align}
where $d(s,z)$ is the linear incomplete gamma function defined by
\begin{equation}
d(s,z)=\int_{0}^{z}e^{-t}t^{s-1}dt=\int_{0}^{\infty}e^{-t}t^{s-1}dt-\int_{z}^{\infty}e^{-t}t^{s-12}dt. \label{24}
\end{equation}
We summarize our results in \eqref{22} and \eqref{23} in the next proposition.
\begin{proposition}
For $n \ge 0$, we have the following identities:
\begin{displaymath}
\mathrm{Bel}_{n,\lambda}^{(p)} = \sum_{k=0}^{n}\sum_{m=0}^{p-1}\frac{(m+1)\binom{p}{m+1}(-1)^{m}S_{2,\lambda}(n,k)}{k+m+1},
\end{displaymath}
and
\begin{displaymath}
\sum_{n=0}^{\infty}\mathrm{Bel}_{n,\lambda}^{(p)}\frac{t^{n}}{n!}=
p\frac{e^{e_{\lambda}(t)-1}}{\big(e_{\lambda}(t)-1\big)^{p}}d\big(p,e_{\lambda}(t)-1\big),
\end{displaymath}
where $d\big(p,e_{\lambda}(t)-1\big)$ is the linear incomplete gamma function in \eqref{24}.
\end{proposition}

From \eqref{13}, we note that
\begin{align}
&\sum_{n=0}^{\infty}\mathrm{Bel}_{n,\lambda}^{(p)}(x)\frac{t^{n}}{n!}=\frac{p!}{x^{p}\big(e_{\lambda}(t)-1\big)^{p}}\sum_{n=p}^{\infty}\frac{x^{n}(e_{\lambda}(t)-1)^{n}}{n!}\label{25} \\
&=\ \frac{p!}{x^{p}(e_{\lambda}(t)-1)^{p}}\sum_{n=0}^{\infty}\frac{x^{n}(e_{\lambda}(t)-1)^{n}}{n!}-\sum_{k=1}^{p}\frac{(p)_{k}}{x^{k}(e_{\lambda}(t)-1)^{k}}\nonumber\\
&=\ \frac{p!}{x^{p}t^{p}}\bigg(\frac{t}{e_{\lambda}(t)-1}\bigg)^{p}\sum_{n=0}^{\infty}\frac{x^{n}(e_{\lambda}(t)-1)^{n}}{n!}-\sum_{k=1}^{p}k!\frac{\binom{p}{k}}{x^{k}}\frac{1}{t^{k}}\bigg(\frac{t}{e_{\lambda}(t)-1}\bigg)^{k}\nonumber\\
&=\ \frac{p!}{x^{p}t^{p}}\sum_{l=0}^{\infty}\beta_{l,\lambda}^{(p)}\frac{t^{l}}{l!}\sum_{m=0}^{\infty}\mathrm{Bel}_{m,\lambda}(x)\frac{t^{m}}{m!}-\sum_{k=1}^{p}k!x^{-k}\binom{p}{k}\frac{1}{t^{k}}\sum_{n=0}^{\infty}\beta_{n,\lambda}^{(k)}\frac{t^{n}}{n!}\nonumber \\
&=\ \frac{p!}{x^{p}}\sum_{n=0}^{\infty}\sum_{m=0}^{n}\binom{n}{m}\mathrm{Bel}_{m,\lambda}(x)\beta_{n-m,\lambda}^{(p)}\frac{t^{n-p}}{n!}-\sum_{k=1}^{p}\sum_{n=0}^{\infty}k!x^{-k}\binom{p}{k}\beta_{n,\lambda}^{(k)}\frac{t^{n-k}}{n!} \nonumber\\
&=\ \frac{p!}{x^{p}}\sum_{n=p}^{\infty}\sum_{m=0}^{n}\binom{n}{m}\mathrm{Bel}_{m,\lambda}(x)\beta_{n-m,\lambda}^{(p)}\frac{t^{n-p}}{n!}-\sum_{k=1}^{p}\sum_{n=k}^{\infty}k!x^{-k}\binom{p}{k}\beta_{n,\lambda}^{(k)}\frac{t^{n-k}}{n!} \nonumber\\
&=\ \sum_{n=0}^{\infty}\sum_{m=0}^{n+p}\binom{n+p}{m}\frac{p!n!}{(n+p)!}\beta_{n+p-m,\lambda}^{(p)}\mathrm{Bel}_{m,\lambda}(x)x^{-p}\frac{t^{n}}{n!} \nonumber\\
&\quad-\sum_{n=0}^{\infty}\sum_{k=1}^{p}\binom{p}{k}\frac{k!n!}{(n+k)!}\beta_{n+k,\lambda}^{(k)}x^{-k}\frac{t^{n}}{n!}\nonumber\\
&=\ \sum_{n=0}^{\infty}\bigg(\sum_{m=0}^{n+p}\frac{\binom{n+p}{m}}{\binom{n+p}{n}}\beta_{n+p-m,\lambda}^{(p)}\mathrm{Bel}_{m,\lambda}(x)x^{-p}-\sum_{k=1}^{p}\frac{\binom{p}{k}}{\binom{n+k}{n}}\beta_{n+k,\lambda}^{(k)}x^{-k}\bigg)\frac{t^{n}}{n!}.\nonumber
\end{align}
Therefore, by comparing the coefficients on both sides of \eqref{25}, we obtain the following theorem.
\begin{theorem}
For $n,p\ge 0$, we have
\begin{displaymath}
x^{p}\mathrm{Bel}_{n,\lambda}^{(p)}(x)=\sum_{m=0}^{n+p}\frac{\binom{n+p}{m}}{\binom{n+p}{n}}\beta_{n+p-m,\lambda}^{(p)}\mathrm{Bel}_{m,\lambda}(x)-\sum_{k=1}^{p}\frac{\binom{p}{k}}{\binom{n+k}{n}}\beta_{n+k,\lambda}^{(p)}x^{p-k}.
\end{displaymath}
In the case of $x=1$, we obtain
\begin{displaymath}
\mathrm{Bel}_{n,\lambda}^{(p)}=\sum_{m=0}^{n+p}\frac{\binom{n+p}{m}}{\binom{n+p}{n}}\beta_{n+p-m,\lambda}^{(p)}\mathrm{Bel}_{m,\lambda}-\sum_{k=1}^{p}\frac{\binom{p}{k}}{\binom{n+k}{n}}\beta_{n+k,\lambda}^{(p)}.
	\end{displaymath}
\end{theorem}
We observe that
\begin{align}
&\bigg(e_{\lambda}^{\lambda-1}(t)\frac{d}{dt}\bigg)^{p}\sum_{n=0}^{\infty}\frac{1}{(n+1)!}\big(1-e_{\lambda}(t)\big)^{n}\label{26}\\
&=(-1)\bigg(e_{\lambda}^{\lambda-1}(t)\frac{d}{dt}\bigg)^{p-1}\sum_{n=1}^{\infty}\frac{1}{n+1}\frac{(1-e_{\lambda}(t))^{n-1}}{(n-1)!}\nonumber \\
&=\ (-1)^{2}\bigg(e_{\lambda}^{\lambda-1}(t)\frac{d}{dt}\bigg)^{p-2}\sum_{n=2}^{\infty}\frac{1}{(n+1)(n-2)!}\big(1-e_{\lambda}(t)\big)^{n-2}\nonumber \\
&=\ \cdots  \nonumber \\
&=\ (-1)^{p}\sum_{n=p}^{\infty}\frac{1}{(n+1)(n-p)!}\big(1-e_{\lambda}(t)\big)^{n-p}\nonumber \\
&=\ (-1)^{p}\sum_{n=0}^{\infty}\frac{1}{(n+p+1)n!}\big(1-e_{\lambda}(t)\big)^{n}. \nonumber
\end{align}
From \eqref{13} and noting that $\frac{1}{(p+n)\binom{p+n-1}{n}}=\frac{(p-1)!n!}{(p+n)!}=\sum_{l=0}^{n}\frac{\binom{n}{l}(-1)^{l}}{p+l}$, we have
\begin{align}
\sum_{n=0}^{\infty}\mathrm{Bel}_{n,\lambda}^{(p)}\frac{t^{n}}{n!}\ &=\ p!\sum_{n=0}^{\infty}\frac{(e_{\lambda}(t)-1)^{n}}{(n+p)!}\ =\ p\sum_{n=0}^{\infty}\frac{(p-1)!n!}{(p+n)!}\frac{1}{n!}\big(e_{\lambda}(t)-1\big)^{n} \label{27}\\
&=\ p\sum_{n=0}^{\infty}\bigg(\sum_{l=0}^{n}\frac{\binom{n}{l}(-1)^{l}}{p+l}\bigg)\frac{1}{n!}\big(e_{\lambda}(t)-1\big)^{n} \nonumber \\
&=\ p\sum_{m=0}^{\infty}\frac{1}{m!}\big(e_{\lambda}(t)-1	\big)^{m}\sum_{l=0}^{\infty}\frac{1}{p+l}\frac{(-1)^{l}}{l!}\big(e_{\lambda}(t)-1\big)^{l} \nonumber \\
&=\ pe^{e_{\lambda}(t)-1}\sum_{l=0}^{\infty}\frac{1}{p+l}\frac{1}{l!}\big(1-e_{\lambda}(t)\big)^{l}.\nonumber
\end{align}
By \eqref{26} and \eqref{27}, we get
\begin{align}
	\sum_{n=0}^{\infty}\mathrm{Bel}_{n,\lambda}^{(p)}\frac{t^{n}}{n!}\ &=\ p!\sum_{n=0}^{\infty}\frac{(e_{\lambda}(t)-1)^{n}}{(n+p)!}\ =\ pe^{e_{\lambda}(t)-1}\sum_{n=0}^{\infty}\frac{1}{p+n}\frac{1}{n!}\big(1-e_{\lambda}(t)\big)^{n} \label{28} \\
	&=\ (-1)^{p-1}pe^{e_{\lambda}(t)-1}\bigg(e_{\lambda}^{\lambda-1}(t)\frac{d}{dt}\bigg)^{p-1}\bigg(\frac{1-e^{1-e_{\lambda}(t)}}{e_{\lambda}(t)-1}\bigg).\nonumber
\end{align}
Therefore, by \eqref{28}, we obtain the following theorem.
\begin{theorem}
	For $p\ge 1$, we have
	\begin{displaymath}
		\sum_{n=0}^{\infty}\mathrm{Bel}_{n,\lambda}^{(p)}\frac{t^{n}}{n!}=(-1)^{p-1}pe^{e_{\lambda}(t)-1}\bigg(e_{\lambda}^{\lambda-1}(t)\frac{d}{dt}\bigg)^{p-1}\bigg(\frac{1-e^{1-e_{\lambda}(t)}}{e_{\lambda}(t)-1}\bigg).
	\end{displaymath}
\end{theorem}
From \eqref{27}, we have
\begin{align}
\sum_{n=0}^{\infty}\mathrm{Bel}_{n,\lambda}^{(p)}\frac{t^{n}}{n!}\ &=\ pe^{e_{\lambda}(t)-1}\sum_{l=0}^{\infty}\frac{1}{p+l}\frac{1}{l!}\big(1-e_{\lambda}(t)\big)^{l} \label{29} \\
&=\ pe^{e_{\lambda}(t)-1}\sum_{l=0}^{\infty}\frac{(-1)^{l}}{p+l}\sum_{m=l}^{\infty}S_{2,\lambda}(m,l)\frac{t^{m}}{m!}\nonumber \\
&=\ p\sum_{k=0}^{\infty}\mathrm{Bel}_{k,\lambda}\frac{t^{k}}{k!}\sum_{m=0}^{\infty} \sum_{l=0}^{m}\frac{(-1)^{l}}{p+l}S_{2,\lambda}(m,l)\frac{t^{m}}{m!} \nonumber\\
&=\ \sum_{n=0}^{\infty}\bigg(p\sum_{m=0}^{n}\sum_{l=0}^{m}\binom{n}{m}\frac{(-1)^{l}}{p+l}S_{2,\lambda}(m,l)\mathrm{Bel}_{n-m,\lambda}\bigg)\frac{t^{n}}{n!}.\nonumber
\end{align}
Therefore, by comparing the coefficients on both sides of \eqref{29}, we obtain the following theorem.
\begin{theorem}
For $n\ge 0$ and $p\ge 1$, we have
\begin{displaymath}
\mathrm{Bel}_{n,\lambda}^{(p)}=p\sum_{m=0}^{n}\sum_{l=0}^{m}\binom{n}{m}\frac{(-1)^{l}}{p+l}S_{2,\lambda}(m,l)\mathrm{Bel}_{n-m,\lambda}.
\end{displaymath}
\end{theorem}
From \eqref{6}, we can derive the following equation.
\begin{align}
\int_{0}^{2\pi}\frac{1}{k!}\big(e_{\lambda}(e^{i\theta})-1\big)^{k}\sin n\theta d\theta & \ =\ \sum_{m=k}^{\infty}S_{2,\lambda}(m,k)\frac{1}{m!}\int_{0}^{2\pi}e^{mi\theta}\sin n\theta d\theta\label{30} \\
&=\ \sum_{m=k}^{\infty}S_{2,\lambda}(m,k)\frac{1}{m!}\int_{0}^{2\pi}\big(\cos m\theta+i\sin m\theta)\sin n\theta d\theta\nonumber \\
&=\ i\sum_{m=k}^{\infty}S_{2,\lambda}(m,k)\frac{1}{m!}\int_{0}^{2\pi}\sin m\theta\sin n\theta d\theta\nonumber\\
&=\ i\frac{S_{2,\lambda}(n,k)}{n!}\int_{0}^{2\pi}\sin^{2}n\theta d\theta\ =\ \frac{i\pi}{n!}S_{2,\lambda}(n,k),\nonumber	
\end{align}
where $i=\sqrt{-1}$, and $n, k$ are integers with $ n>0, \,0\le k\le n$. \par
Thus, by \eqref{30}, we obtain the following lemma.
\begin{lemma}
For $n,k\in\mathbb{Z}$ with $n>0,\,0\le k\le n$, we have
\begin{displaymath}
S_{2,\lambda}(n,k)=\frac{n!}{\pi}\mathrm{Im}\int_{0}^{2\pi}\frac{1}{k!}\big(e_{\lambda}(e^{i\theta})-1\big)^{k}\sin n\theta d\theta.
\end{displaymath}
\end{lemma}
By Lemma 9, we get
\begin{align}
\int_{0}^{2\pi}e^{e_{\lambda}(e^{i\theta})-1}\sin n\theta d\theta&=\sum_{k=0}^{\infty}\frac{1}{k!}\int_{0}^{2 \pi}\big(e_{\lambda}(e^{i\theta})-1\big)^{k}\sin n\theta d\theta \label{31} \\
&=\ \sum_{k=0}^{\infty}\sum_{m=k}^{\infty}S_{2,\lambda}(m,k)\frac{1}{m!}\int_{0}^{2\pi}e^{mi\theta}\sin n\theta d\theta \nonumber \\
&=\ \sum_{m=0}^{\infty}\sum_{k=0}^{m}\frac{S_{2,\lambda}(m,k)}{m!}\int_{0}^{2\pi}\big(\cos m\theta+i\sin m\theta\big)\sin n\theta d\theta \nonumber\\
&=\ i\sum_{m=0}^{\infty}\sum_{k=0}^{m}\frac{S_{2,\lambda}(m,k)}{m!}\int_{0}^{2\pi}\sin m\theta\sin n\theta d\theta \nonumber \\
&=\ i\frac{\pi}{n!}\sum_{k=0}^{n}S_{2,\lambda}(n,k)\ =\ i\frac{\pi}{n!}\mathrm{Bel}_{n,\lambda},\nonumber
\end{align}
where $n$ is a positive integer. \par
Therefore, by \eqref{31}, we obtain the following corollary.
\begin{corollary}
For $n > 0$, we have
\begin{displaymath}
\frac{n!}{\pi}\mathrm{Im}\int_{0}^{2\pi}e^{e_{\lambda}(e^{i\theta})-1}\sin n\theta d\theta=\mathrm{Bel}_{n,\lambda}.
\end{displaymath}
\end{corollary}
From \eqref{13}, we note that
\begin{align}
&p!\int_{0}^{2\pi}\bigg(\frac{e^{e_{\lambda}(e^{i\theta})-1}}{(e_{\lambda}(e^{i \theta})-1)^{p}}-\sum_{l=0}^{p-1}\frac{(e_{\lambda}(e^{i\theta})-1)^{l-p}}{l!}\bigg)\sin n\theta d\theta \label{32}\\
&\quad=\ \sum_{m=0}^{\infty}\mathrm{Bel}_{m,\lambda}^{(p)}\frac{1}{m!}\int_{0}^{2\pi}e^{im\theta}\sin n\theta d\theta\nonumber\\
&\quad=\ \sum_{m=0}^{\infty}\mathrm{Bel}_{m,\lambda}^{(p)}\frac{1}{m!}\int_{0}^{2\pi}\big(\cos m\theta+i\sin m\theta\big)\sin n\theta d\theta \nonumber \\
&\quad=\ i\sum_{m=0}^{\infty}\mathrm{Bel}_{m,\lambda}^{(p)}\frac{1}{m!}\int_{0}^{2\pi}\sin m\theta\sin n\theta d\theta\ =\ i\frac{\pi}{n!}\mathrm{Bel}_{n,\lambda}^{(p)}, \nonumber
\end{align}
where $n$ is a positive integer. \par
Therefore, by \eqref{32}, we obtain the following theorem.
\begin{theorem}
For $n > 0$, we have
\begin{displaymath}
\mathrm{Bel}_{n,\lambda}^{(p)}=\frac{n!p!}{\pi}\mathrm{Im}\int_{0}^{2\pi}\bigg(\frac{e^{e_{\lambda}(e^{i\theta})-1}}{(e_{\lambda}(e^{i\theta})-1)^{p}}-\sum_{l=0}^{p-1}\frac{\big(e_{\lambda}(e^{i\theta})-1\big)^{l-p}}{l!}\bigg)\sin n\theta d\theta.
\end{displaymath}
\end{theorem}
From \eqref{13} and \eqref{15}, we have
\begin{align}
e_{\lambda}^{\lambda-1}(t)\frac{d}{dt}\sum_{n=0}^{\infty}\mathrm{Bel}_{n,\lambda}^{(p)}(1)\frac{t^{n}}{n!}&=e_{\lambda}^{\lambda-1}(t)\frac{d}{dt}p!\sum_{n=0}^{\infty}\frac{(e_{\lambda}(t)-1)^{n}}{(n+p)!}\label{33} \\	
&=\ p!\sum_{n=0}^{\infty}\frac{n+1}{(n+p+1)!}\big(e_{\lambda}(t)-1\big)^{n} \nonumber\\
&=\ p\sum_{n=0}^{\infty}B(n+2,p)\frac{1}{n!}\big(e_{\lambda}(t)-1\big)^{n}\nonumber\\
&=\ p\int_{0}^{1}x\sum_{n=0}^{\infty}\frac{1}{n!}x^{n}\big(e_{\lambda}(t)-1\big)^{n}(1-x)^{p-1}dx \nonumber\\
&=\ p\int_{0}^{1}e^{x(e_{\lambda}(t)-1)}x(1-x)^{p-1}dx.\nonumber
\end{align}
It is easy to see that
\begin{align}
e_{\lambda}^{\lambda-1}(t)\frac{d}{dt}\sum_{n=0}^{\infty}\mathrm{Bel}_{n,\lambda}^{(p)}\frac{t^{n}}{n!}\ &=\ \sum_{m=0}^{\infty}\mathrm{Bel}_{m+1,\lambda}^{(p)}\frac{t^{m}}{m!}\sum_{l=0}^{\infty}(\lambda-1)_{l,\lambda}\frac{t^{l}}{l!}\label{34}\\
&=\ \sum_{n=0}^{\infty}\bigg(\sum_{m=0}^{n}\binom{n}{m}\mathrm{Bel}_{m+1,\lambda}^{(p)}(\lambda-1)_{n-m,\lambda}\bigg)\frac{t^{n}}{n!}. \nonumber
\end{align}
On the other hand, we also have
\begin{align}
p\int_{0}^{1}e^{x(e_{\lambda}(t)-1)}x(1-x)^{p-1}dx
&=-p\int_{0}^{1}e^{x(e_{\lambda}(t)-1)}(1-x)^p dx+p\int_{0}^{1}e^{x(e_{\lambda}(t)-1)}(1-x)^{p-1}dx \label{35}	\\
&=-\frac{p}{p+1}(p+1)!\sum_{n=0}^{\infty}\frac{(e_{\lambda}(t)-1)^{n}}{(n+p+1)!}+p!\sum_{n=0}^{\infty}\frac{(e_{\lambda}(t)-1)^{n}}{(n+p)!}\nonumber\\
&=\sum_{n=0}^{\infty}\bigg(-\frac{p}{p+1}\mathrm{Bel}_{n,\lambda}^{(p+1)}+\mathrm{Bel}_{n,\lambda}^{(p)}\bigg)\frac{t^{n}}{n!}.\nonumber
\end{align}
Therefore, by \eqref{33}, \eqref{34} and \eqref{35}, we obtain the following theorem.
\begin{theorem}
For $p\ge 0$ and $n \ge 2$, we have
\begin{displaymath}
\mathrm{Bel}_{n+1,\lambda}^{(p)}=(n+1-n\lambda)\mathrm{Bel}_{n,\lambda}^{(p)}-\frac{p}{p+1}\mathrm{Bel}_{n,\lambda}^{(p)}-\sum_{m=0}^{n-2}\binom{n}{m}\mathrm{Bel}_{m+1,\lambda}^{(p)}(\lambda-1)_{n-m,\lambda}.
\end{displaymath}
In the case of $p=0$, we obtain
\begin{align*}
\mathrm{Bel}_{n+1,\lambda}&=(n+1-n\lambda)\mathrm{Bel}_{n,\lambda}-\sum_{m=0}^{n-2}\binom{n}{m}\mathrm{Bel}_{m+1,\lambda}(\lambda-1)_{n-m,\lambda} \\
&=(n+1-n\lambda)\mathrm{Bel}_{n,\lambda}-\sum_{m=0}^{n-1}\binom{n}{m-1}\mathrm{Bel}_{m,\lambda}(\lambda-1)_{n-m+1,\lambda}.
\end{align*}
\end{theorem}
It is known that the degenerate Stirling polynomials of the second kind are defined by
\begin{equation}
(y+x)_{n,\lambda}=\sum_{k=0}^{n}S_{2,\lambda}(n,k|x)(y)_{k},\quad(n\ge 0),\quad(\mathrm{see}\ [10,13]).\label{36}	
\end{equation}
Thus, we note that
\begin{equation}
\sum_{n=0}^{\infty}\bigg(\frac{1}{n!}\big(e_{\lambda}(t)-1\big)^{n}e_{\lambda}^{x}(t)\bigg)(y)_{n}\ =\ e_{\lambda}^{y+x}(t)\ =\ \sum_{k=0}^{\infty}(y+x)_{k,\lambda}\frac{t^{k}}{k!} \label{37}
\end{equation}
\begin{displaymath}
=\sum_{k=0}^{\infty}\sum_{n=0}^{k}S_{2,\lambda}(k,n|x)(y)_{n}\frac{t^{k}}{k!}=\sum_{n=0}^{\infty}\bigg(\sum_{k=n}^{\infty}S_{2,\lambda}(k,n|x)\frac{t^{k}}{k!}\bigg)(y)_{n}.
\end{displaymath}
By comparing the coefficients on both sides of \eqref{37}, we get
\begin{equation}
\frac{1}{n!}\big(e_{\lambda}(t)-1\big)^{n}e_{\lambda}^{x}(t)=\sum_{k=n}^{\infty}S_{2,\lambda}(k,n|x)\frac{t^{k}}{k!},\quad(n\ge 0). \label{38}	
\end{equation}
From \eqref{5}, we note that
\begin{align}
\sum_{n=0}^{\infty}(x+y)_{n,\lambda}\frac{t^{n}}{n!}\ &=\ e_{\lambda}^{x+y}(t)\ =\ \sum_{i=0}^{\infty}(x)_{i,\lambda}\frac{t^{i}}{i!}\sum_{l=0}^{\infty}(y)_{l,\lambda}\frac{t^{l}}{l!} \label{39} \\
&=\ \sum_{n=0}^{\infty}\sum_{i=0}^{n}\binom{n}{i}(x)_{n-i,\lambda}(y)_{i,\lambda}\frac{t^{n}}{n!}.\nonumber
\end{align}
By \eqref{36} and \eqref{39}, we get
\begin{align}
\sum_{l=0}^{n}S_{2,\lambda}(n,l|x)(y)_{l}\ &=\ (y+x)_{n,\lambda}\ =\ \sum_{i=0}^{n}\binom{n}{i}(x)_{n-i,\lambda}(y)_{i,\lambda}\nonumber \\
&=\ \sum_{i=0}^{n}\binom{n}{i}(x)_{n-i,\lambda}\sum_{l=0}^{i}S_{2,\lambda}(i,l)(y)_{l} \label{40}\\
&=\ \sum_{l=0}^{n}\bigg(\sum_{i=l}^{n}\binom{n}{i}(x)_{n-i,\lambda}S_{2,\lambda}(i,l)\bigg)(y)_{l}.\nonumber
\end{align}
By comparing the coefficients on both sides of \eqref{40}, we obtain the following theorem.
\begin{theorem}
For $n,l\ge 0$, we have
\begin{displaymath}
S_{2,\lambda}(n,l|x)=\sum_{i=l}^{n}\binom{n}{i}(x)_{n-i,\lambda}S_{2,\lambda}(i,l).
\end{displaymath}
\end{theorem}
For $p\ge 0$, we define {\it{truncated degenerate modified Bell polynomials}} $B_{n,\lambda}^{(p)}(x)$ by
\begin{equation}
\frac{p!}{\big(e_{\lambda}(t)-1\big)^{p}}\bigg(e^{e_{\lambda}(t)-1}-\sum_{l=0}^{p-1}\frac{(e_{\lambda}(t)-1)^{l}}{l!}\bigg)e_{\lambda}^{x}(t)=\sum_{n=0}^{\infty}B_{n,\lambda}^{(p)}(x)\frac{t^{n}}{n!}. \label{41}	
\end{equation}
For $x=1$, $B_{n,\lambda}^{(p)}=B_{n,\lambda}^{(p)}(1)$ are called {\it{truncated degenerate modified Bell numbers}}.
From \eqref{41}, we note that
\begin{align}
\sum_{n=0}^{\infty}B_{n,\lambda}^{(p)}(x)\frac{t^{n}}{n!}\ &=\ p!\sum_{n=0}^{\infty}\frac{(e_{\lambda}(t)-1)^{n}}{(n+p)!}e_{\lambda}^{x}(t) \label{42} \\
&=\ \sum_{m=0}^{\infty}\mathrm{Bel}_{n,\lambda}^{(p)}\frac{t^{m}}{m!}\sum_{l=0}^{\infty}(x)_{l,\lambda}\frac{t^{l}}{l!}\nonumber\\
&=\ \sum_{n=0}^{\infty}\bigg(\sum_{m=0}^{n}\binom{n}{m}\mathrm{Bel}_{n,\lambda}^{(p)}(x)_{n-m,\lambda}\bigg)\frac{t^{n}}{n!}.\nonumber	
\end{align}
Thus we obtain
\begin{displaymath}
B_{n,\lambda}^{(p)}(x)=\sum_{m=0}^{n}\binom{n}{m}\mathrm{Bel}_{n,\lambda}^{(p)}(x)_{n-m,\lambda}.
\end{displaymath}
Again, from \eqref{41} we observe that
\begin{align}
\sum_{n=0}^{\infty}B_{n,\lambda}^{(p)}(x)\frac{t^n}{n!}\
&=\ \sum_{k=0}^{\infty}\frac{p!k!}{(k+p)!}\frac{1}{k!}\big(e_{\lambda}(t)-1\big)^{k}e_{\lambda}^{x}(t) \label{43} \\
&=\ \sum_{k=0}^{\infty}\frac{1}{\binom{k+p}{p}}\sum_{n=k}^{\infty}S_{2,\lambda}(n,k|x)\frac{t^{n}}{n!} \nonumber \\
&=\ \sum_{n=0}^{\infty}\bigg(\sum_{k=0}^{n}\frac{S_{2,\lambda}(n,k|x)}{\binom{k+p}{p}}\bigg)\frac{t^{n}}{n!}. \nonumber
\end{align}
Therefore, by \eqref{43}, we obtain the following theorem.
\begin{theorem}
For $n,p\ge 0$, we have
\begin{displaymath}
B_{n,\lambda}^{(p)}(x)=\sum_{k=0}^{n}\frac{1}{\binom{k+p}{p}}S_{2,\lambda}(n,k|x).
\end{displaymath}
\end{theorem}
From \eqref{41}, we have
\begin{align}
p\int_{0}^{1}(1-y)^{p-1}e_{\lambda}^{x}(t)e^{y(e_{\lambda}(t)-1)}dy\ &=\ pe_{\lambda}^{x}(t)\sum_{n=0}^{\infty}\frac{(e_{\lambda}(t)-1)^{n}}{n!}\int_{0}^{1}(1-y)^{p-1}y^{n}dy
\label{44}	\\
&=\ pe_{\lambda}^{x}(t)\sum_{n=0}^{\infty}\frac{(e_{\lambda}(t)-1)^{n}}{n!}B(n+1,p) \nonumber \\
&=\ p!\sum_{n=0}^{\infty}\frac{(e_{\lambda}(t)-1)^{n}}{(n+p)!}e_{\lambda}^{x}(t) \nonumber \\
&=\ \sum_{n=0}^{\infty}B_{n,\lambda}^{(p)}(x)\frac{t^{n}}{n!}. \nonumber
\end{align}
Thus, we get
\begin{equation}
p\int_{0}^{1}(1-y)^{p-1}e_{\lambda}^{x}(t)e^{y(e_{\lambda}(t)-1}dy=\sum_{n=0}^{\infty}B_{n,\lambda}^{(p)}(x)\frac{t^{n}}{n!}.\label{45}	
\end{equation}
From \eqref{45}, we can easily derive the following equation
\begin{equation}
B_{n,\lambda}^{(p)}(x)=p\sum_{m=0}^{n}\binom{n}{m}(x)_{n-m,\lambda}\int_{0}^{1}(1-y)^{p-1}\mathrm{Bel}_{m,\lambda}(y)dy,\quad(n\ge 0), \label{46}	
\end{equation}
By \eqref{9}, we easily get
\begin{equation}
\mathrm{Bel}_{m,\lambda}(x)=e^{-x}\sum_{k=0}^{\infty}\frac{(k)_{m,\lambda}}{k!}x^k.\label{47}
\end{equation}
Thus, by \eqref{46} and \eqref{47}, we get
\begin{align}
B_{n,\lambda}^{(p)}(x)\ &=\ p\sum_{m=0}^{n}\binom{n}{m}(x)_{n-m,\lambda}\int_{0}^{1}(1-y)^{p-1}\mathrm{Bel}_{m,\lambda}(y)dy \label{48}\\
&=\ p\sum_{k=0}^{\infty}\frac{(x+k)_{n,\lambda}}{k!}\sum_{m=0}^{\infty}\frac{(-1)^{m}}{m!}\int_{0}^{1}y^{m+k}(1-y)^{p-1}dy \nonumber \\
&=\ \sum_{k=0}^{\infty}\sum_{m=0}^{\infty}\frac{(x+k)_{n,\lambda}}{k!}\frac{(-1)^{m}}{m!\binom{m+k+p}{p}}.\nonumber
\end{align}
On the other hand,
\begin{align}
&\ p\sum_{m=0}^{n}\binom{n}{m}(x)_{n-m,\lambda}\int_{0}^{1}(1-y)^{p-1}\mathrm{Bel}_{m,\lambda}(y)dy \label{49}\\
&=\ p\sum_{m=0}^{n}\binom{n}{m}(x)_{n-m,\lambda}\sum_{k=0}^{m}S_{2,\lambda}(m,k)\int_{0}^{1}(1-y)^{p-1}y^{k}dy  \nonumber\\
&=\ \sum_{m=0}^{n}\sum_{k=0}^{m}\binom{n}{m}(x)_{n-m,\lambda}S_{2,\lambda}(m,k)\binom{p+k}{k}^{-1}.\nonumber
\end{align}
Therefore, by \eqref{48} and \eqref{49}, we obtain the following theorem.
\begin{theorem}
For $n,p\ge 0$, we have
\begin{displaymath}
\sum_{k=0}^{\infty}\sum_{m=0}^{\infty}\frac{(x+k)_{n,\lambda}}{k!}\frac{(-1)^{m}}{m!\binom{m+k+p}{p}}=\sum_{m=0}^{n}\sum_{k=0}^{m}\frac{\binom{n}{m}}{\binom{p+k}{k}}S_{2,\lambda}(m,k)(x)_{n-m,\lambda}.
	\end{displaymath}
\end{theorem}
From \eqref{48}, we note that
\begin{equation}
B_{n,\lambda}^{(p)}(x)\ =\ p\sum_{k=0}^{\infty}\frac{(x+k)_{n,\lambda}}{k!}\int_{0}^{1}e^{-y}y^{k}(1-y)^{p-1}dy,\quad (n\ge 0). \label{50}
\end{equation}
Thus, we get
\begin{align}
\ B_{n+1,\lambda}^{(p)}(x)\ &=\ p\sum_{k=0}^{\infty}\frac{(x+k)_{n+1,\lambda}}{k!}\int_{0}^{1}e^{-y}y^{k}(1-y)^{p-1}dy \label{51}\\
&\ =\ (x-n\lambda)p\sum_{k=0}^{\infty}\frac{(x+k)_{n,\lambda}}{k!}\int_{0}^{1}e^{-y}y^{k}(1-y)^{p-1}dy\nonumber \\
&\quad +p\sum_{k=1}^{\infty}\frac{(x+k)_{n,\lambda}}{(k-1)!}\int_{0}^{1}e^{-y}y^{k}(1-y)^{p-1}dy \nonumber \\
&=\ (x-n\lambda)\mathrm{Bel}_{n,\lambda}^{(p)}(x)+p\sum_{k=0}^{\infty}\frac{(x+k+1)_{n,\lambda}}{k!}\int_{0}^{1}y^{k}ye^{-y}(1-y)^{p-1}dy \nonumber \\
&=\ (x-n\lambda)B_{n,\lambda}^{(p)}(x)-\sum_{j=0}^{n}\binom{n}{j}(1)_{n-j,\lambda}p\sum_{k=0}^{\infty}\frac{(x+k)_{j,\lambda}}{k!}\int_{0}^{1}(1-y)^{p}e^{-y}y^{k}dy\nonumber \\
&\quad +\sum_{j=0}^{n}\binom{n}{j}(1)_{n-j,\lambda}p\sum_{k=0}^{\infty}\frac{(x+k)_{j,\lambda}}{k!}\int_{0}^{1}(1-y)^{p-1}e^{-y}y^{k}dy \nonumber \\
&=\ (x-n\lambda)B_{n,\lambda}^{(p)}(x)-\sum_{j=0}^{n}\binom{n}{j}(1)_{n-j,\lambda}\bigg(\frac{p}{p+1}B_{j,\lambda}^{(p+1)}(x)-B_{j,\lambda}^{(p)}(x)\bigg).\nonumber
\end{align}
Therefore, by \eqref{51}, we obtain the following theorem.
\begin{theorem}
For $n,p\ge 0$, we have
\begin{displaymath}
B_{n+1,\lambda}^{(p)}(x)=(x-n\lambda)B_{n,\lambda}^{(p)}(x)-\sum_{j=0}^{n}\binom{n}{j}(1)_{n-j,\lambda}\big(\frac{p}{p+1}B_{j,\lambda}^{(p+1)}(x)-B_{j,\lambda}^{(p)}(x)\big).
\end{displaymath}
\end{theorem}

\section{Further Remark}
It is well known that $X$ is the beta random variable with parameters $\alpha>0$ and $\beta>0$ if the probability density function of $X$ is given by
\begin{displaymath}
f(x)=\left\{\begin{array}{ccc}
\frac{1}{B(\alpha,\beta)}x^{\alpha-1}(1-x)^{\beta-1}, & \textrm{if $0\le x\le 1$} \\
0, &  \textrm{otherwise.}
\end{array}\right.
\end{displaymath}
The beta random variable $X$ with parameters $\alpha>0$ and $\beta>0$ is denoted by $X\sim\mathrm{Beta}(\alpha,\beta)$, (see [22]). \par
Let $g(x)$ be a real valued function, and let $f(x)$ be the probability density function of the continuous random variable $X$. Then the expectation of $g(X)$ is defined by
\begin{displaymath}
	E\big[g(X)\big]=\int_{-\infty}^{\infty}g(x)f(x)dx,\quad(\mathrm{see}\ [21]).
\end{displaymath}
For $X\sim\mathrm{Beta}(\alpha,\beta)$, we have
\begin{align*}
	E\big[e^{tX}\big]\ &=\ \int_{-\infty}^{\infty}f(x)e^{tx}\ = \
 \frac{1}{B(\alpha,\beta)}\int_{0}^{1}x^{\alpha-1}(1-x)^{\beta-1}e^{tx}dx\\
 &=\ \frac{1}{B(\alpha,\beta)}\sum_{n=0}^{\infty}\int_{0}^{1}x^{n+\alpha-1}(1-x)^{\beta-1}dx\frac{t^{n}}{n!} \\
&=\ \sum_{n=0}^{\infty}\frac{B(n+\alpha,\beta)}{B(\alpha,\beta)}\frac{t^{n}}{n!}.
\end{align*}
Thus, we get
\begin{displaymath}
E\big[X^{n}\big]\ =\ \frac{B(n+\alpha,\beta)}{B(\alpha,\beta)}=\frac{\binom{n+\alpha-1}{n}}{\binom{n+\alpha+\beta-1}{n}},\quad (n\ge 0).
\end{displaymath}
For $X\sim\mathrm{Beta}(1,p)$, we have
\begin{align}
&\ E\big[e^{X(e_{\lambda}(t)-1)}\big]\ =\ p\int_{0}^{1}(1-x)^{p-1}e^{x(e_{\lambda}(t)-1)}dx\nonumber\\
&\ =\ \sum_{n=0}^{\infty}\frac{(e_{\lambda}(t)-1)^{n}}{n!}p\int_{0}^{1}(1-x)^{p-1}x^{n}dx\ =\ p!\sum_{n=0}^{\infty}\frac{(e_{\lambda}(t)-1)^{n}}{(n+p)!}\label{52}\\
&\ =\ \sum_{n=0}^{\infty}\mathrm{Bel}_{n,\lambda}^{(p)}\frac{t^{n}}{n!}.\nonumber 	
\end{align}
On the other hand,
\begin{align}
E\big[e^{X(e_{\lambda}(t)-1)}\big]\ &=\ \sum_{k=0}^{\infty}E\big[X^{k}\big]\frac{1}{k!}\big(e_{\lambda}(t)-1\big)^{k}\label{53} \\
&=\ \sum_{k=0}^{\infty}E[X^{k}]\sum_{n=k}^{\infty}S_{2,\lambda}(n,k)\frac{t^{n}}{n!} \nonumber \\
&=\ \sum_{n=0}^{\infty}\bigg(\sum_{k=0}^{n}S_{2,\lambda}(n,k)E[X^{k}]\bigg)\frac{t^{n}}{n!}.\nonumber
\end{align}
From \eqref{52} and \eqref{53}, we have
\begin{displaymath}
\sum_{k=0}^{n}E[X^{k}]S_{2,\lambda}(n,k)=\mathrm{Bel}_{n,\lambda}^{(p)},\quad(n\ge 0),
\end{displaymath}
where $X\sim\mathrm{Beta}(1,p)$.

\section{Conclusion}

Truncated polynomials have been shown to play an important role in various areas.
In this paper, we introduced truncated degenerate Bell polynomials and numbers and investigated their various properties. In more detail, for the truncated degenerate Bell polynomials and numbers we obtained explicit expressions, identities involving other special polynomials, integral representations, Dobinski-like formula and  expressions of the generating function in terms of differential operators and linear incomplete gamma function. In addition, we introduced  truncated degenerate modified Bell polynomials and numbers and got similar results for those polynomials. As an application of our results, we showed that the truncated degenerate Bell numbers can be expressed as a finite sum involving moments of the beta random variable with some parameters. \par
We would like to conclude this section by displaying the different expressions for the truncated degenerate Bell numbers that we have obtained in this paper.

\begin{align*}
\mathrm{Bel}_{n,\lambda}^{(p)} &=p\int_{0}^{1}\mathrm{Bel}_{n,\lambda}(x)(1-x)^{p-1}dx \\
&=\sum_{k=0}^{\infty}\sum_{l=0}^{\infty}(-1)^{l}\frac{\binom{k+l}{l}}{\binom{k+l+p}{p}}\frac{(k)_{n,\lambda}}{(k+l)!}\\
&= \sum_{k=0}^{n}\sum_{m=0}^{p-1}\frac{(m+1)\binom{p}{m+1}(-1)^{m}S_{2,\lambda}(n,k)}{k+m+1}\\
&=p\sum_{m=0}^{n}\sum_{l=0}^{m}\binom{n}{m}\frac{(-1)^{l}}{p+l}S_{2,\lambda}(m,l)\mathrm{Bel}_{n-m,\lambda}\\
&=\frac{n!p!}{\pi}\mathrm{Im}\int_{0}^{2\pi}\bigg(\frac{e^{e_{\lambda}(e^{i\theta})-1}}{(e_{\lambda}(e^{i\theta})-1)^{p}}-\sum_{l=0}^{p-1}\frac{\big(e_{\lambda}(e^{i\theta})-1\big)^{l-p}}{l!}\bigg)\sin n\theta d\theta \\
&=\sum_{k=0}^{n}E[X^{k}]S_{2,\lambda}(n,k) ,\quad(n\ge 0),
\end{align*}
where $X\sim\mathrm{Beta}(1,p)$.

\bigskip

\noindent{\bf{Acknowledgements}}\\
The authors would like to thank
Jangjeon Research Institute for Mathematical Sciences for the support of this research.\\

\bigskip

\noindent{\bf{Funding}}\\
Not applicable. \\

\bigskip

\noindent{\bf{Availability of data and materials}}\\
Not applicable. \\

\bigskip

\noindent{\bf{Ethics approval and consent to participate}}\\
All authors reveal that there is no ethical problem in the publishing this paper.

\bigskip

\noindent{\bf{Competing interests}}\\
The authors declare that they have no competing interests.

\bigskip

\noindent{\bf{Consent for publication}}\\
All authors want to publish this paper in this journal.

\bigskip

%\noindent{\bf{Authors’ contributions}}\\
%TK and DSK conceived of the framework and structured the whole paper; TK and DSK wrote the paper; L-CJ and HL checked the results of the paper and typed the paper; DSK, TK and L-CJ completed the revision of the article. All %authors have read and agreed with the published version of the manuscript. \\

\bigskip

%\noindent{\bf{Author details}}


\begin{thebibliography}{9}
\bibitem{1}
Brillhart, J. \emph{Mathematical Notes: Note on the single variable Bell polynomials,} Amer. Math. Monthly  \textbf{74} (1967), no. 6, 695–-696.
\bibitem{2}
Carlitz, L. \emph{Single variable Bell polynomials,} Collect. Math. \textbf{14} (1962), 13-–25.
\bibitem{3}
Carlitz, L. \emph{Arithmetic properties of the Bell polynomials,} J. Math. Anal. Appl. \textbf{15} (1966), 33-–52.
\bibitem{4}
Charalambides, Ch. A. \emph{Minimum variance unbiased estimation for the zero class truncated bivariate Poisson and logarithmic series distributions,} Metrika  \textbf{31} (1984), no. 2, 115-–123.
\bibitem{5}
Cohen, A.; Clifford, Jr. \emph{Estimation in the truncated Poisson distribution when zeros and some ones are missing,} J. Amer. Statist. Assoc. \textbf{55} (1960), 342–-348.
\bibitem{6}
Comtet, L. \emph{Advanced combinatorics. The art of finite and infinite expansions,} Revised and enlarged edition. D. Reidel Publishing Co., Dordrecht, 1974.
\bibitem{7}
Dattoli, G.; Cesarano, C.; Sacchetti, D. \emph{A note on truncated polynomials,} Appl. Math. Comput. \textbf{134} (2003), no. 2-3, 595–-605.
\bibitem{8}
Guettai, G.; Laissaoui, D.; Rahmani, M.; Sebaoui, M. \emph{On poly-Bell numbers and polynomials,} arXiv:1812.04136 .
\bibitem{9}
Kim, H. K.;  Baek, H.; Lee, D. S. \emph{A note on truncated degenerate exponential polynomials,} Research Gate(preprint: December 2020). https://www.researchgate.net/publication/346581812
\bibitem{10}
Kim, D. S.; Kim, T. \emph{A note on a new type of degenerate Bernoulli numbers,} Russ. J. Math. Phys. \textbf{27} (2020), no. 2, 227-–235.
\bibitem{11}
Kim, D. S.; Kim, T. \emph{A note on polyexponential and unipoly functions.} Russ. J. Math. Phys. \textbf{26} (2019), no. 1, 40-–49.
\bibitem{ }
Kim, D. S.; Kim, T. An umbral calculus approach to Bernoulli-Pad\'e polynomials, Mathematical analysis, approximation theory and their applications, 363-–382, Springer Optim. Appl., 111, Springer, [Cham], 2016.
\bibitem{12}
Kim, T. \emph{A note on degenerate Stirling polynomials of the second kind,} Proc. Jangjeon Math. Soc. \textbf{20} (2017), no. 3, 319-–331.
\bibitem{13}
Kim, T.; Kim, D. S. \emph{Degenerate zero-truncated Poisson random variables,} arXiv:1911.13227.
\bibitem{14}
Kim, T.; Kim, D. S. \emph{Degenerate polyexponential functions and degenerate Bell polynomials,} J. Math. Anal. Appl. \textbf{487} (2020), no. 2, 124017, 15 pp.
\bibitem{15}
Kim, T.; Kim, D. S. \emph{Note on the degenerate gamma function,} Russ. J. Math. Phys. \textbf{27} (2020), no. 3, 352-–358.
\bibitem{16}
Kim, T.; Kim, D. S. \emph{Degenerate Laplace transform and degenerate gamma function,} Russ. J. Math. Phys. \textbf{24} (2017), no. 2, 241-–248.
\bibitem{17}
Kim, T.; Kim, D. S.; Dolgy, D. V. \emph{On partially degenerate Bell numbers and polynomials,} Proc. Jangjeon Math. Soc. \textbf{20} (2017), no. 3, 337-–345.
\bibitem{18}
Kim, T.; Kim, D. S.; Jang, L.-C.; Kwon, H.-I. \emph{Extended degenerate Stirling numbers of the second kind and extended degenerate Bell polynomials,} Util. Math. \textbf{106} (2018), 11-–21.
\bibitem{19}
Roman, S. \emph{The umbral calculus,} Pure and Applied Mathematics, 111. Academic Press, Inc. [Harcourt Brace Jovanovich, Publishers], New York, 1984.
\bibitem{20}
Ross, S. M. \emph{Introduction to probability models,} Twelfth edition of [MR0328973], Academic Press, London, 2019.
\bibitem{21}
Springer, M. D.; Thompson, W. E. \emph{The distribution of products of beta, gamma and Gaussian random variables,} SIAM J. Appl. Math. \textbf{18} (1970), 721-–737.
\bibitem{22}
Srivastava, H.M.; Araci, S.; Khan,W.A.; Acikgoz, M. \emph{A note on the truncated-exponential based Apostol-type polynomials,} Symmetry \textbf{11} (2019), no. 4, 538, 15pp.
\bibitem{23}
Whittaker, E. T.; Watson, G. N. \emph{A course of modern analysis, An introduction to the general theory of infinite processes and of analytic functions: with an account of the principal transcendental functions,} Fourth edition, Reprinted Cambridge University Press, New York, 1962.
\bibitem{24}
Yasmin, G.; Khan, S.; Ahmad, N. \emph{Operational methods and truncated exponential-based Mittag-Leffler
polynomials,} Mediterr. J. Math. \textbf{13} (2016), 1555-1569.
\end{thebibliography}
\end{document}